\documentclass[12pt, reqno]{amsart}
\usepackage{amsmath, amsthm, amscd, amsfonts, amssymb, graphicx, color}
\usepackage[bookmarksnumbered, colorlinks, plainpages]{hyperref}

\newtheorem{theorem}{Theorem}[section]
\newtheorem{lemma}[theorem]{Lemma}

\newtheorem{corollary}[theorem]{Corollary}
\theoremstyle{definition}
\newtheorem{definition}[theorem]{Definition}
\newtheorem{example}[theorem]{Example}

\theoremstyle{remark}
\newtheorem{remark}[theorem]{Remark}
\numberwithin{equation}{section}

\begin{document}

\title[Some fixed point theorems on partial cone metric spaces]
{Some fixed point theorems for multi-valued mappings on partial cone metric spaces}

\author[T. L. Shateri]{T. L. Shateri}
\address{Tayebe Lal Shateri \\ Department of Mathematics and Computer
Sciences, Hakim Sabzevari University, Sabzevar, P.O. Box 397, IRAN}
\email{ \rm t.shateri@hsu.ac.ir; shateri@ualberta.ca}
\author[H. Isik]{H. Isik}
\address{Huseyin Isik \\ Faculty of Science, University of Gazi, 06500-Teknikokullar, Ankara, Turkey.}
\email{ \rm isikhuseyin76@gmail.com}
\thanks{*The corresponding author:
isikhuseyin76@gmail.com (Huseyin Isik)}
 \maketitle

\begin{abstract}
In this paper, we study the fixed point theory for multi-valued mappings on partial cone metric spaces. We prove an analogous to the well-known Kannan$'s$ fixed point theorem and Chatterjea$'s$ fixed point theorem for multi-valued mappings on partial cone metric spaces.
\vskip 3mm
\end{abstract}

\section{Introduction and preliminaries}\vskip 2mm
In 1994, Matthews \cite{MAT} introduced the notion of partial metric space as a part of the study of denotational semantics of data ﬂow networks in computer science. In 2009, Bukatin established the precise relationship between partial metric spaces and the so called weightable quasi-metric space \cite{BU}.  
The concept of a partial metric, and any concept related to a partial metric play a very important role not only in pure mathematics but also in other branches of science involving mathematics especially in computer science, information science, and biological science.
In partial metric spaces, the self-distance for any point need not be equal to zero.\\
In 1980, Rzepecki \cite{RZ} introduced a generalized metric $d_E$ on a set $\mathcal X$ in a way that $d_E:\mathcal X\times \mathcal X\to P$, replacing the set of real numbers with a Banach space $E$ where $P$ is a normal cone in $E$ with a partial order $\leq$. Lin \cite{LI} considered the notion of cone metric spaces by replacing real numbers with a cone $P$ in which it is called a $K$-metric. Many authors to investigate topological properties of cone metric spaces (see \cite{CA,TU}).

Fixed point theorems are the basic mathematical tools used in showing the
existence of solution concept in such diverse ﬁelds as biology, chemistry, economics, engineering, and game theory. It is well-known that the Banach contraction principle is a fundamental result in ﬁxed point theory, which has been used, and extended in many diﬀerent directions. Mattews proved the Banach ﬁxed point theorem \cite{BA} in  Partial metric space. Some fixed point theorems of contractive mappings for partial metric spaces and cone metric spaces have proved by many authors (see \cite{AB1,AB2,HU}). Complete partial metric spaces constitute a suitable framework to model several distinguished examples of the theory of computation and also to model metric spaces via domain theory (see \cite{CI,RO,SC}). Partial cone metric space have investigated by Mahlotra, et al. \cite{MAH}. They proved some fixed point theorems in this space.\\
The study of fixed points for multi-valued contractions using the Hausdorff metric was initiated by Nadler \cite{NA}. The theory of multi-valued maps has application in control theory, convex optimization, differential equations and economics (see \cite{DA}).

In this paper, by using some ideas of \cite{AY,DE} we prove an analogous to the  Kannan fixed point theorem \cite{KA} and Chatterjea's fixed point theorem \cite{CH} for multi-valued mappings on partial cone metric spaces.

The following definitions and results will be needed in the sequel.
\begin{definition}(Partial metric space) A partial metric on a non-empty set $\mathcal X$ is a function $p:\mathcal X\times \mathcal X \to \mathbb R^+$ such that for all $x,y,z\in \mathcal X$ the followings hold\\
$(p_1)\; 0\leq p(x,x)\leq p(x,y),$\\
$(p_2)\; x=y \quad \text{if and only if}\quad p(x,x)=p(x,y)=p(y,y),$\\
$(p_3)\; p(x,y)=p(y,x),$\\
$(p_4)\; p(x,y)\leq p(x,z)+p(z,y)-p(z,z).$\\
Then the pair $(\mathcal X,p)$ is called a partial metric space. 
\end{definition}
It is clear that if $p(x,y)=0$, then $(p1)$ and $(p2)$ imply that $x=y$. But if $x=y$, $p(x,y)$ may not be $0$. A basic example of a partial metric space is the pair $(\mathbb R^+,p)$, where $p(x,y)=\max{x,y}$ for all $x,y\in \mathbb R^+$.\\
Let $E$ be a real Banach space and $P$ a subset of $E$. $P$ is called a cone if it satisfies the followings\\
$(C_1)\; P \;\text{is closed, non-empty and}\; P\neq \{0\},$\\
$(C_2)\; ax+by\in P \quad \text{for all}\; x,y\in P\; \text{and non-negative real numbers}\; a,b,$\\
$(C_3)\; P\cap (-P)=\{0\}.$\\
For a given cone $P\subseteq E$, we can deﬁne a partial ordering $\leq $ on $E$ with respect to $P$ by $x\leq y$ if and only if $y-x\in P$. We write $x<y$ to 
indicate that $x\leq y$ but $x\neq y$, while $x\ll y$ will stand for $y-x\in  intP$, in which $intP$ denotes the interior of $P$. The cone $P$ is called normal if there is a constant number $M>0$ such that for all $x,y\in E$ where $0\leq x\leq y$ implies $\|x\|\leq M\|y\|$. The least positive number satisfying above is called the normal constant of $P$. 

Let $E$ be a Banach space, $P$ a cone in $E$ with $intP\neq\Phi$ and $\leq$ is partial ordering with respect to $P$. 
\begin{definition}(Cone metric space) Let $\mathcal X$ be a non-empty set. The  mapping $d:\mathcal X\times\mathcal X\to E$ is said to be a cone metric on $\mathcal X$ if for all $x,y,z\in X$ the followings hold\\
$(CM_1)\; 0\leq d(x,y)\;\text{and}\; d(x,y)=0\; \text{if and only if}\; x=y,$\\
$(CM_2)\; d(x,y)=d(y,x),$\\
$(CM_3)\; d(x,y)\leq d(x,z)+d(z,y).$\\
Then $(\mathcal X,d)$ is called a cone metric space.
\end{definition} 
Mahlotra and et al. \cite{MAH} and S$\ddot{o}$nmez \cite{SO} introduced the notion of partial cone metric space and its topological characterization. They also developed some ﬁxed point theorems in this generalized setting. We now state the deﬁnition of partial cone metric space.
\begin{definition} (Partial cone metric space) A partial cone metric on a non-empty set $\mathcal X$ is a function $p:\mathcal X\times\mathcal X\to\mathbb E$ such that for all $x,y,z\in X$\\ 
$(PCM_1)\; 0\leq p(x,x)\leq p(x,y),$\\
$(PCM_2)\; x=y\quad \text{if and only if}\quad p(x,x)=p(x,y)=p(y,y),$\\
$(PCM_3)\; p(x,y)=p(y,x),$\\
$(PCM_4)\; p(x,y)\leq p(x,z)+p(z,y)-p(z,z).$
\end{definition}
A partial cone metric space is a pair $(\mathcal X,p)$ such that $\mathcal X$ is a non-empty set and $p$ is a partial cone metric on $\mathcal X$. It is clear that, if $p(x,y)=0$, then $(PCM_1)$ and $(PCM_2)$ imply that $x=y$. But the converse is not true ingeneral. A cone metric space is a partial cone metric space, but there exists partial cone metric spaces which are not cone metric space. we give the following example from \cite{SO}.
\begin{example}
Let $E=\mathbb R^2$, $P=\{(x,y)\in E:\; x,y\geq 0\}$ and $\mathcal X=\mathbb R^+$ and $p:\mathcal X\times \mathcal X\to E$ deﬁned by $p(x,y)=(\max\{x,y\},\alpha\max\{x,y\})$ where $\alpha\geq 0$ is a constant. Then $(\mathcal X,p)$ is a partial cone metric space which is not a cone metric space.
\end{example}
Following, We give some properties of partial cone metric spaces, for more details see \cite{SO}.
\begin{theorem}\cite[Theorem 4]{SO} 
 Every partial cone metric space $(\mathcal X,p)$ is a topological space.
\end{theorem}
\begin{definition}
Let $(\mathcal X,p)$ be a partial cone metric space. Let $\{x_n\}$ be a sequence in $\mathcal X$ and $x\in \mathcal X$.\\
$(i)$ $\{x_n\}$ is said to be convergent to $x$ and $x$ is called a limit of $\{x_n\}$ if
$$\lim_{n\to\infty}p(x_n,x)=\lim_{n\to\infty}p(x_n,x_n)=p(x,x).$$
$(ii)$ Let $(\mathcal X,p)$ be a partial cone metric space. $\{x_n\}$ be a sequence in$\mathcal  X$. $\{x_n\}$ is Cauchy sequence if there is $x\in P$ such that for every $\varepsilon >0$ there is $N$ such that for all $n,m> N$, $\|p(x_n,x_m)-x\|<\varepsilon$.\\
$(i)\; (\mathcal X,p)$ is said to be complete if every Cauchy sequence in $(\mathcal X,p)$ is convergent in $(\mathcal X,p)$.
\end{definition}
\begin{theorem}\label{eq2}\cite[Theorem 6]{SO} 
Let $(\mathcal X,p)$ be a partial cone metric space, $A\subset \mathcal X$ and $a\in \mathcal X$. Then $a\in \bar{A}$ if and only if $p(a,A)=p(a,a)$.
\end{theorem}
\begin{corollary}\cite[Corollary 3]{SO} 
Every closed subset of a partial cone metric space is complete.
\end{corollary}
Suppose $(\mathcal X,p)$ is a partial cone metric space, then
$$d(x,y)=2p(x,y)-p(x,x)-p(y,y)\quad (x,y\in\mathcal X),$$
is a cone metric on $\mathcal X$. The following result is given.
\begin{theorem}\label{th0}\cite[Theorem 5]{SO}
Let $(\mathcal X,p)$ be a partial cone metric space and let $P$ be a normal cone. If $\{x_n\}$ is a Cauchy sequence in $(\mathcal X,p)$, then it is a Cauchy sequence in the cone metric space $(\mathcal X,d)$.
\end{theorem}
In fact, the last theorem, establish a correspondence between a partial cone metric space and a cone metric space.
\section{\bf Main results}
Let $(\mathcal X,p)$ be a partial cone metric space. Throughout this section, we assume that $CB^p(\mathcal X)$ is a family of all nonempty, closed and bounded subsets of the partial cone metric space $(\mathcal X,p)$, induced by the cone partial metric $p$.\\
We now assert the following deﬁnitions from \cite{AY}.
\begin{definition}
Let $\mathcal X$ be a nonempty set and $2^{\mathcal X}$ denote the set of all nonempty subsets of $\mathcal X$. An element $x\in \mathcal X$ is said to be a fixed point of a multi-valued mapping $T:\mathcal X\to 2^{\mathcal X}$ if $x\in Tx$.
\end{definition}
For $A,B\in CB^p(\mathcal X)$ and $x\in \mathcal X$, define
\begin{align*}
&p(x,A)=\inf\{p(x,a):a\in A\},\quad \delta_p(A,B)=\sup\{p(a,B):a\in A\},\quad \\&\delta_p(B,A)=\sup\{p(b,A):b\in B\}\text{and}\quad H_p(A,B)=\max\{\delta_p(A,B),\delta_p(B,A)\}.
\end{align*}
 The proof of the following lemma is similar to \cite[Lemma 3.1]{AY}. 
\begin{lemma}\label{lem1}
Let $(\mathcal X,p)$ be a partial cone metric space, $A,B\in CB^p(\mathcal X)$ and $h>1$. For $a\in A$, there exists $b=b(a)\in B$ such that $$p(a,b)\leq hH_p(A,B).$$
\end{lemma}
Now, we state and prove our main results.\\
In the following theorem we prove  Kannan's fixed point theorem \cite{KA} for multi-valued mappings on complete partial cone metric spaces. 
\begin{theorem}\label{th1}
Let $(\mathcal X,p)$ be a complete partial cone metric space, $P$ be a normal
cone with normal constant $M$. Let $T:\mathcal X\to CB^p(\mathcal X)$ be a  multi-valued mapping such that
\begin{equation}\label{eq11}
H_p(Tx,Ty)\leq \lambda\left[p(Tx,x)+p(Ty,y)\right],
\end{equation}
for all $x,y \in \mathcal X$, where $\lambda\in(0,\frac{1}{2})$, is a constant. Then $T$ has a ﬁxed point in $\mathcal X$
\begin{proof}
Let $x_0\in \mathcal X$ and $x_1\in Tx_0$. From Lemma \eqref{lem1} with $h=\frac{1}{2\lambda}$, there exists  $x_2\in Tx_1$ such that  
\begin{align*}
p(x_2,x_1)&\leq hH_p(Tx_1,Tx_0)\\&\leq h\lambda\left[p(Tx_1,x_1)+p(Tx_0,x_0)\right]\\
&\leq h\lambda\left[p(x_2,x_1)+p(x_1,x_0)\right],
\end{align*}
and so $p(x_2,x_1)\leq \frac{h\lambda}{1-h\lambda}p(x_1,x_0)$. For $x_2\in Tx_1$, there exists $x_3\in Tx_2$ such that $p(x_3,x_2)\leq \frac{h\lambda}{1-h\lambda}p(x_2,x_1)$. Set $k=\frac{h\lambda}{1-h\lambda}$, since $\lambda\in(0,\frac{1}{2})$, hence $k<1$ and $p(x_3,x_2)\leq k^2p(x_1,x_0)$. Continuing this process, we obtain a sequence $\{x_n\}$ in $\mathcal X$ such that $$x_{n+1}\in Tx_n \quad and \quad p(x_{n+1},x_n)\leq k p(x_n,x_{n-1})\leq k^np(x_1,x_0).$$ We show that $\{x_n\}$ is a Cauchy sequence. For $m>n$ we have
\begin{align*}
p(x_m,x_n)&\leq p(x_m,x_{m-1})+p(x_{m-1},x_{m-2})+\ldots+p(x_{n+1},x_n)-\sum_{i=1}^{m-n-1}p(x_{m-i},x_{m-i})\\
&\leq p(x_m,x_{m-1})+p(x_{m-1},x_{m-2})+\ldots+p(x_{n+1},x_n)\\
&\leq(k^{m-1}+k^{m-2}+\ldots k^n)p(x_1,x_0)\\
&=\frac{k^n}{1-k}p(x_1,x_0)
\end{align*}
and so
$$\|p(x_m,x_n)\|\leq \frac{k^n}{1-k}M\|p(x_1,x_0)\|.$$
Hence $p(x_m,x_n)\to 0$ as $m,n\to \infty$. Therefore $\{x_n\}$ is a Cauchy sequence, and by completeness of $\mathcal X$, there exists $x\in \mathcal X$ such that $x_n\to x$ as $n\to \infty$. Hence
\begin{equation}\label{eq12}
p(x,x)=\lim_{n\to\infty}p(x_n,x)=\lim_{n\to\infty}p(x_n,x_n)=0.
\end{equation}
By \eqref{eq11} we obtain
\begin{align*}
p(x,Tx)&\leq p(x_n,Tx)+p(x_n,x)-p(x_n,x_n)\\
&\leq \delta_p(Tx_{n-1},Tx)+p(x_n,x)\\
&\leq H_p(Tx_{n-1},Tx)+p(x_n,x)\\
&\leq \lambda\left[p(Tx_{n-1},x_{n-1})+p(Tx,x)\right]+p(x_n,x)\\
&\leq \lambda\left[p(x_n,x_{n-1})+p(Tx,x)\right]+p(x_n,x),
\end{align*}
therefore
$$(1-\lambda)p(x,Tx)\leq\lambda p(x_n,x_{n-1})+p(x_n,x),$$
gives
$$p(x,Tx)\leq \frac{1}{1-\lambda}\left[\lambda p(x_n,x_{n-1})+p(x_n,x)\right],$$
so
$$\|p(x,Tx)\|\leq \frac{M}{1-\lambda}\left[\lambda \|p(x_n,x_{n-1})\|+\|p(x_n,x)\|\right]\to 0.$$
Hence $\|p(x,Tx)\|=0$, and so $p(x,Tx)=0$. Therefore, from \eqref{eq12} we get
$p(x,x)=p(x,Tx)$, and Theorem \ref{eq2} implies that $x\in\overline{Tx}=Tx$. This completes the proof.
\end{proof}
\end{theorem}
By using \cite[Example 3.3]{AY}, we give the following very simple illustrative example.
\begin{example}
Let $E=\mathbb R^2$, $P=\{(x,y)\in E:\; x,y\geq 0\}$ and let $\mathcal X=\{0,1,4\}$ be endowed with the partial cone metric $$p(x,y)=\left(\frac{1}{4}|x-y|,\frac{1}{2}\max\{x,y\}\right)\quad (x,y\in\mathcal X).$$ 
Note that $p$ is not a cone metric, because $p(1,1)=(0,\frac{1}{4})$. 
Since $d(x,y)=(\frac{1}{2}|x-y|,\frac{1}{2}|x-y|)$, Theorem \ref{th0} implies that $(\mathcal X,p)$ is a complete partial cone metric space. It is clear that $\{0\}$ and $\{0,1\}$ are bounded subsets in $(\mathcal X,p)$. Moreover, we have
\begin{align*}
x\in\overline{\{0\}}&\Leftrightarrow\quad p(x,\{0\})=p(x,x)\\
&\Leftrightarrow\quad(\frac{1}{4}x,\frac{1}{2}x)=(0,\frac{1}{2}x)\\
&\Leftrightarrow\quad x=0\\
&\Leftrightarrow\quad x\in\{0\}.\\
\text{Hence}\; \{0\}\; \text{is closed in}\; (\mathcal X,p).\; \text{Also}\\
x\in\overline{\{0,1\}}&\Leftrightarrow\quad p(x,\{0,1\})=p(x,x)\\
&\Leftrightarrow\quad\min\left\{(\frac{1}{4}x,\frac{1}{2}x),\left(\frac{1}{4}|x-1|,\frac{1}{2}\max\{x,1\}\right)\right\}=(0,\frac{1}{2}x)\\
&\Leftrightarrow\quad x\in\{0,1\}.
\end{align*}
Therefore $\{0,1\}$ is closed in $(\mathcal X,p)$.

Now, consider the mapping $T:\mathcal X\to CB^p(\mathcal X)$ defined as
$$T(0)=T(1)=\{0\}\quad \text{and}\quad T(4)=\{0,1\}.$$
We prove that the condition \ref{eq11} is satisfied with $\lambda=\frac{1}{3}$. We have the following cases:

$(i)$ if $x,y\in\{0,1\}$, then
$$H_p(T(x),T(y))=H_p(\{0\},\{0\})=(0,0).$$
$(ii)$ if $x\in\{0,1\}$ and $y=4$, we have
\begin{align*}
H_p(T(0),T(4))&=H_p(T(1),T(4))\\
&=H_p(\{0\},\{0,1\})\\
&=\max\left\{p(0,\{0,1\}),\max\{p(0.0),p(1,0)\}\right\}\\
&=\max\left\{p(0,0),p(1,0)\right\}\\
&=(\frac{1}{4},\frac{1}{2}),
\end{align*}
on the other hand 
$$p(T(0),0)+p(T(4),4)=(\frac{3}{4},2)\quad \text{and}\quad p(T(1),1)+p(T(4),4)=(1,\frac{3}{2}).$$
Therefore we get
$$H_p(T(x),T(4))\leq\lambda\left[p(T(x),x)+p(T(4),4)\right],$$
for $x\in\{0,1\}$.
$(iii)$ if $x=y=4$, then
\begin{align*}
H_p(T(4),T(4))&=H_p(\{0,1\},\{0,1\})\\
&=\max\left\{p(x,\{0,1\})x\in\{0,1\}\right\}\\
&=\max\left\{p(0,0),p(1,1)\right\}\\
&=(0,\frac{1}{2})\\
&\leq\frac{1}{3}(\frac{3}{2},4)\\
&=\lambda\left[p(T(4),4)+p(T(4),4)\right].
\end{align*}
Therefore $T$ satisfies in the all hypothesis of Theorem \ref{th1}, and $x=0$ is a fixed point of $T$.
\end{example}
We now prove the ﬁxed point theorem due to Chatterjea \cite{CH} for multi-valued mappings on partial cone metric spaces.
\begin{theorem}\label{th2}
Let $(\mathcal X,p)$ be a complete partial cone metric space, $P$ be a normal
cone with normal constant $M$. Let $T:\mathcal X\to CB^p(\mathcal X)$ be a  multi-valued mapping such that
\begin{equation}\label{eq21}
H_p(Tx,Ty)\leq \lambda\left[p(Tx,y)+p(Ty,x)\right],
\end{equation}
for all $x,y \in \mathcal X$, where $\lambda\in(0,\frac{1}{2})$, is a constant. Then $T$ has a ﬁxed point in $\mathcal X$
\begin{proof}
Let $x_0\in \mathcal X$ and $x_1\in Tx_0$. From Lemma \eqref{lem1} with $h=\frac{1}{2\lambda}$, there exists  $x_2\in Tx_1$ such that  
\begin{align*}
p(x_2,x_1)&\leq hH_p(Tx_1,Tx_0)\leq h\lambda\left[p(Tx_1,x_0)+p(Tx_0,x_1)\right]\\
&\leq h\lambda\left[p(x_2,x_0)+p(x_1,x_1)\right]\\
&\leq h\lambda\left[p(x_2,x_1)+p(x_1,x_0)-p(x_1,x_1)+p(x_1,x_1)\right]\\
&= h\lambda\left[p(x_2,x_1)+p(x_1,x_0)\right],
\end{align*}
and so $p(x_2,x_1)\leq \frac{h\lambda}{1-h\lambda}p(x_1,x_0)$. For $x_2\in Tx_1$, there exists $x_3\in Tx_2$ such that $p(x_3,x_2)\leq \frac{h\lambda}{1-h\lambda}p(x_2,x_1)$. Set $k=\frac{h\lambda}{1-h\lambda}$, since $\lambda\in(0,\frac{1}{2})$, hence $k<1$ and $p(x_3,x_2)\leq k^2p(x_1,x_0)$. Continuing this process, we obtain a sequence $\{x_n\}$ in $\mathcal X$ such that $$x_{n+1}\in Tx_n \quad and \quad p(x_{n+1},x_n)\leq k p(x_n,x_{n-1})\leq k^np(x_1,x_0).$$ 
Proceeding as the proof of Theorem \ref{th1} we see that $\{x_n\}$ is a Cauchy sequence. By completeness of $\mathcal X$, there exists $x\in \mathcal X$ such that $x_n\to x$ as $n\to \infty$. Hence
\begin{equation}\label{eq22}
p(x,x)=\lim_{n\to\infty}p(x_n,x)=\lim_{n\to\infty}p(x_n,x_n)=0.
\end{equation}
By \eqref{eq21} we obtain
\begin{align*}
p(x,Tx)&\leq p(x_n,Tx)+p(x_n,x)-p(x_n,x_n)\\
&\leq \delta_p(Tx_{n-1},Tx)+p(x_n,x)\\
&\leq H_p(Tx_{n-1},Tx)+p(x_n,x)\\
&\leq \lambda\left[p(Tx_{n-1},x)+p(Tx,x_{n-1})\right]+p(x_n,x)\\
&\leq \lambda\left[p(x_n,x)+p(Tx,x_{n-1})\right]+p(x_n,x)\\
&\leq \lambda\left[p(x_n,x)+p(Tx,x)+p(x,x_{n-1}-p(x,x))\right]+p(x_n,x)\\
&\leq \lambda\left[p(x_n,x)+p(Tx,x)+p(x,x_{n-1})\right]+p(x_n,x),
\end{align*}
therefore
$$(1-\lambda)p(x,Tx)\leq\lambda\left[p(x_n,x)+p(x,x_{n-1})\right]+p(x_n,x),$$
gives
$$p(x,Tx)\leq \frac{1}{1-\lambda}\left[\lambda p(x_n,x)+\lambda p(x,x_{n-1})+p(x_n,x)\right],$$
so
$$\|p(x,Tx)\|\leq \frac{M}{1-\lambda}\left[\lambda \|p(x_n,x)\|+\lambda \|p(x,x_{n-1})\|+\|p(x_n,x)\|\right]\to 0.$$
Hence $\|p(x,Tx)\|=0$, and so $p(x,Tx)=0$. Therefore, from \eqref{eq22} we get
$p(x,x)=p(x,Tx)$, and Theorem \ref{eq2} implies that $x\in\overline{Tx}=Tx$. This completes the proof.
\end{proof}
\end{theorem}
\begin{example}
Let $E=\mathbb R^2$, $P=\{(x,y)\in E:\; x,y\geq 0\}$ and let $\mathcal X=\{0,1,2\}$ be endowed with the partial cone metric $p:\mathcal X\times\mathcal X\to\mathbb E$ defined by
\begin{align*}
&p(0,0)=p(1,1)=(0,0), \quad p(2,2)=(\frac{1}{4},0),\\
&p(0,1)=p(1,0)=(\frac{1}{6},0),\\
&p(0,2)=p(2,0)=(\frac{7}{10},0),\\
&p(1,2)=p(2,1)=(\frac{1}{2},0).\\
\end{align*}
Define the mapping $T:\mathcal X\to CB^p(\mathcal X)$ defined as
$$T(0)=T(1)=\{0\}\quad \text{and}\quad T(2)=\{0,1\}.$$
We prove that the condition \ref{eq11} is satisfied with $\lambda=\frac{1}{4}$. We have the following cases:

$(i)$ if $x,y\in\{0,1\}$, then
$$H_p(T(x),T(y))=H_p(\{0\},\{0\})=(0,0).$$
$(ii)$ if $x\in\{0,1\}$ and $y=2$, we obtain
\begin{align*}
H_p(T(0),T(2))&=H_p(T(1),T(2))\\
&=H_p(\{0\},\{0,1\})\\
&=\max\left\{p(0,0),p(1,0)\right\}\\
&=(\frac{1}{6},0),
\end{align*}
on the other hand
\begin{align*}
&p(T(0),2)=p(T(1),2)=p(0,2)=(\frac{7}{10},0),\\
&p(T(2),0)=p(\{0,1\},0)=(0,0)\\
\text{and}\\
&p(T(2),1)=p(\{0,1\},1)=p(1,1)=(0,0).
\end{align*}
Therefore we get
$$H_p(T(x),T(2))\leq\lambda\left[p(T(x),2)+p(T(2),x)\right],$$
for $x\in\{0,1\}$.
$(iii)$ if $x=y=2$, then
\begin{align*}
H_p(T(2),T(2))&=H_p(\{0,1\},\{0,1\})\\
&=\max\left\{p(x,\{0,1\})x\in\{0,1\}\right\}\\
&=\max\left\{p(0,0),p(1,1)\right\}\\
&=(0,0)\\
&\leq\frac{1}{4}(1,0)\\
&=\lambda\left[p(T(2),2)+p(2,T(2))\right].
\end{align*}
Hence $T$ satisfies in the all hypothesis of Theorem \ref{th2}, and $x=0$ is a fixed point of $T$.
\end{example}
Following is the Reich’s type contraction mapping \cite{RE} considered here to prove an another ﬁxed point theorem for multi-valued mappings on partial cone metric spaces.
\begin{theorem}\label{th3}
Let $(\mathcal X,p)$ be a complete partial cone metric space, $P$ be a normal
cone with normal constant $M$. Let $T:\mathcal X\to CB^p(\mathcal X)$ be a  multi-valued mapping such that
\begin{equation}\label{eq31}
H_p(Tx,Ty)\leq \alpha p(Tx,x)+\beta p(Ty,y)+\gamma p(x,y),
\end{equation}
for all $x,y \in \mathcal X$, where $0\leq \alpha+\beta+\gamma <1$ and $\alpha,\beta,\gamma$ are non-negative. Then $T$ has a ﬁxed point in $\mathcal X$
\begin{proof}
Let $x_0\in \mathcal X$ and $x_1\in Tx_0$. From Lemma \eqref{lem1} with $h=\frac{1}{\alpha+\beta+\gamma}$, there exists  $x_2\in Tx_1$ such that  
\begin{align*}
p(x_2,x_1)&\leq hH_p(Tx_1,Tx_0)\leq h\left[\alpha p(Tx_1,x_1)+\beta p(Tx_0,x_0)+\gamma p(x_1,x_0)\right]\\
&\leq h\left[\alpha p(x_2,x_1)+\beta p(x_1,x_0)+\gamma p(x_1,x_0)\right],
\end{align*}
and so $p(x_2,x_1)\leq \frac{h(\beta+\gamma)}{1-h\alpha}p(x_1,x_0)$. For $x_2\in Tx_1$, there exists $x_3\in Tx_2$ such that $p(x_3,x_2)\leq \frac{h(\beta+\gamma)}{1-h\alpha}p(x_2,x_1)$. Set $k=\frac{h(\beta+\gamma)}{1-h\alpha}$, since $\alpha+\beta+\gamma <1$, hence $k<1$ and $p(x_3,x_2)\leq k^2p(x_1,x_0)$. Continuing this process, we obtain a sequence $\{x_n\}$ in $\mathcal X$ such that $$x_{n+1}\in Tx_n \quad and \quad p(x_{n+1},x_n)\leq k p(x_n,x_{n-1})\leq k^np(x_1,x_0).$$ Proceeding as the proof of Theorem \ref{th1} we see that $\{x_n\}$ is a Cauchy sequence. By completeness of $\mathcal X$, there exists $x\in \mathcal X$ such that $x_n\to x$ as $n\to \infty$. Hence
\begin{equation}\label{eq32}
p(x,x)=\lim_{n\to\infty}p(x_n,x)=\lim_{n\to\infty}p(x_n,x_n)=0.
\end{equation}
By \eqref{eq31} we obtain
\begin{align*}
p(x,Tx)&\leq p(x_n,Tx)+p(x_n,x)-p(x_n,x_n)\\
&\leq \delta_p(Tx_{n-1},Tx)+p(x_n,x)\\
&\leq H_p(Tx_{n-1},Tx)+p(x_n,x)\\
&\leq \alpha p(x_n,x_{n-1})+\beta p(x,Tx)+\gamma p(x_{n-1},x_{n-1})+p(x_n,x),
\end{align*}
therefore
$$(1-\beta)p(x,Tx)\leq \alpha p(x_n,x_{n-1})+\gamma p(x_{n-1},x_{n-1})+p(x_n,x),$$
gives
$$p(x,Tx)\leq \frac{1}{1-\beta}\left[\alpha p(x_n,x_{n-1})+\gamma p(x_{n-1},x_{n-1})+p(x_n,x)\right],$$
so
$$\|p(x,Tx)\|\leq \frac{1}{1-\beta}\left[\alpha \|p(x_n,x_{n-1})\|+\gamma \|p(x_{n-1},x_{n-1})\|+\|p(x_n,x)\|\right]\to 0.$$
Hence $\|p(x,Tx)\|=0$, and so $p(x,Tx)=0$. Therefore, from \eqref{eq32} we get
$p(x,x)=p(x,Tx)$, and Theorem \ref{eq31} implies that $x\in\overline{Tx}=Tx$. This completes the proof.
\end{proof}
\end{theorem}
\begin{remark}
$(a)$ If in Theorem \ref{th3}, put $\alpha=\beta=\lambda$ and $\gamma=0$, then we get Theorem \ref{th1}.\\
\end{remark}
Putting $\alpha=\beta=0$ and $\gamma=k$ in Theorem \ref{th3}, then we get the following corollary.
\begin{corollary}\cite[Theorem 3.2]{AY}
Let $(\mathcal X,p)$ be a complete partial metric space. If\\ $T:\mathcal X\to CB^p(\mathcal X)$ is a multi-valued mapping such that for all $x,y\in\mathcal X$, we have $$H_p(Tx,Ty)\leq kp(x,y),$$
where $k\in(0,1)$. Then $T$ has a fixed point.
\end{corollary}


\begin{thebibliography}{99}
\bibitem{AB1}  M. Abbas, B. E. Rhoades, \emph{Fixed and periodic point results in cone metric space}, Appl. Math. Lett., \textbf{22}(4) (2009), 511--515.

\bibitem{AB2} M. Abbas, G. Jungck, \emph{Common ﬁxed point results for noncommuting mappings without continuity in cone metric spaces}, J. Math. Anal. Appl., \textbf{341} (2008), 416--420.

\bibitem{AY} H. Aydi, M. Abbas and C. Vetro, \emph{Partial Hausdorff metric and Nadler's fixed point theorem on partial metric spaces}, Top. Appl., 159 (2012), 3234--3242.

\bibitem{BA}  S. Banach, \emph{Sur les operations dans les ensembles abstraits et leur application aux quations intgrales}, Fund. Math., \textbf{3} (1922),133--181 (French).

\bibitem{BU} ] M. Bukatin, R. Kopperman, Steve Matthews, and Homeira Pajoohesh, \emph{Partial Metric Spaces}, 
Amer. Math. Monthly, \textbf{116} (2009), 708--718.


\bibitem{CA}  H. Cakalli, A. S$\ddot{o}$nmez and C.Genc, \emph{On Equivalence of Topological Vector Space Valued Cone Metric Spaces and Metric spaces}, Appl. Math. Lett. \textbf{25} (2012), 429--433.

\bibitem{CH} S.K. Chatterjee, \emph{Fixed point theorems}, Rend. Acad. Bulgare Sc., \textbf{25} (1972), 727--730.

\bibitem{CI} L. $\acute{C}$iri$\acute{c}$, B. Samet, H. Aydi and C. Vetro, \emph{Common fixed points of generalized contractions on partial metric spaces and an application}, Appl. Math. Comput, \textbf{218} (2011), 2398--2406.

\bibitem{DA} B. Damjanovic, B. Samet and C. Vetro, \emph{Common fixed point theorems for multi-valued maps}, Acta. Math. Sci.Ser. B Engl. Ed., \textbf{32} (2012), 818--824.

\bibitem{DE} D. Dey and M. Saha, \emph{Partial cone metric space and some fixed point theorems}, TWMS J. App. Eng. Math., \textbf{3}(1) (2013), 1--9.


\bibitem{HU}  L.-G. Huang, X. Zhang, \emph{Cone metric spaces and ﬁxed point theorems of contractive mappings}, J. Math. Anal. Appl., \textbf{332}(2) (2007), 1468--1476.

\bibitem{KA} R. Kannan, \emph{Some results on ﬁxed points}, Bull. Calcutta Math. Soc., \textbf{60} (1968), 71--76.

\bibitem{LI} S. D. Lin, \emph{A common ﬁxed point theorem in abstract spaces}, Indian J. Pure Appl. Math., \textbf{18}(8) (1987), 685--690.

\bibitem{MAH}  S. K. Mahlotra, S. Shukla, R. Sen, N. Verma, \emph{Fixed point theorems in partial cone metric spaces}, Inter. J. Math. Arch., \textbf{2}(4)  (2011), 610--616.

\bibitem{MAT}  S. G. Matthews, \emph{Partial Metric Topology}, in: Proceedings of the 8th Summer Conference on Topology and its Applications, 728, Annals of The New york Academy of Sciences, (1994), 183--197. MR 98d:54054.

\bibitem{NA} S.B. Nadler, \emph{Multivalued contraction mappings}, Pacific J. Math., \textbf{30} (1969), 475--488.

\bibitem{RE} S. Reich, \emph{Kannan’s ﬁxed point thorem}, Boll. Un. Math. Ital., \textbf{4} (1971), 1--11.

\bibitem{RO} S. Romaguera, \emph{A Kirk type characterization of completeness for partial metric spaces}, Fixed Point Theory Appl., \textbf{2010} (2010), Article ID 493298, 6 pp.

\bibitem{RZ}  B. Rzepecki, \emph{On ﬁxed point theorems of Maia type}, Publ. Inst. Math., \textbf{28}(42) (1980), 179--186.

\bibitem{SC} M.P. Schellekens, \emph{The correspondence between partial metrics and semivaluations}, Theoret. Comput. Sci., \textbf{315} (2004), 135--149. 

\bibitem{SO} A. S$\ddot{o}$nmez, \emph{Fixed point theorems in partial cone metric spaces}, arXiv:1101.2741v1 [math.GN]. 

\bibitem{TU} D. Turkoglu, M. Abuloha, \emph{Cone metric spaces and ﬁxed point theorems in diametrically contractive mappings}, Acta Math. Sin. Engl. Series, \textbf{26}(3) (2010), 489--496.


\end{thebibliography}
\end{document}